\documentclass[a4paper,12pt,english]{article}
\usepackage{babel}
\usepackage[latin1]{inputenc}
\author{Andrea Tellini \\ \small{Departamento de Matem\'atica Aplicada} \\ \small{Universidad Complutense de Madrid} \\ \small{Plaza de Ciencias 3 - 28040 Madrid, Spain} \\ \texttt{andrea.tellini@mat.ucm.es}}
\title{\textbf{Imperfect bifurcations via topological methods in superlinear indefinite problems}\thanks{This research has been supported by Project
MTM2012-30669 and Grant BES-2010-039030 of the
Spanish Ministry of Economy and Competitiveness.}}
\date{\today}
\usepackage{amsmath,amsfonts,amssymb}
\usepackage{yhmath,mathrsfs,yfonts,arcs}
\usepackage{amsthm}
\usepackage{graphicx}

\usepackage{hyperref}
\hypersetup{colorlinks=false,       
    linkcolor=red,          
    citecolor=green,        
    filecolor=magenta,      
    urlcolor=cyan           
    }

\newtheorem{theorem}{Theorem}[section]

\newtheorem{proposition}[theorem]{Proposition}

 \DeclareMathOperator{\sign}{sign}

\newcommand{\field}[1] {\mathbb{#1}}
\newcommand{\N}{\field{N}}

\newcommand{\R}{\field{R}}

\def\a{\alpha}

\def\d{\delta}
\def\g{\gamma}
\def\G{\Gamma}
\def\l{\lambda}

\def\O{\Omega}
\def\p{\partial}

\def\t{\theta}

\def\v{\varphi}
\def\ov{\overline}

\def\ua{\uparrow}
\def\da{\downarrow}

\newcommand{\mc}{\mathcal}

\begin{document}
\maketitle

\begin{abstract}
In \cite{LTZ} the structure of the bifurcation diagrams of a class of superlinear indefinite problems with a symmetric weight was ascertained, showing that they consist of a primary branch and secondary loops bifurcating from it. In \cite{LT} it has been proved that, when the weight is asymmetric, the bifurcation diagrams are no longer connected since parts of the primary branch and loops of the symmetric case form an arbitrarily high number of isolas. In this work we give a deeper insight on this phenomenon, studying how the secondary bifurcations break as the weight is perturbed from the symmetric situation. Our proofs rely on the approach of \cite{LTZ, LT}, i.e. on the construction of certain Poincar\'e maps and the study of how they vary as some of the parameters of the problems change, constructing in this way the bifurcation diagrams.
\end{abstract}



\setcounter{equation}{0}
\section{Introduction}
\label{section1}
This work deals with positive solutions of the boundary value problem
\begin{equation}
\label{eq1}
  \left\{ \begin{array}{l}
  -u''=\l u+a(t)u^p \quad \hbox{in}\;\; \O:=(0,1)\cr
  u(0)=u(1)=M\end{array}\right.
\end{equation}
where $M\in(0,+\infty]$, $p>1$, $\l<0$ are constants and $a(t)$ is a \emph{weight} function indefinite in sign which is defined as
\begin{equation*}
  a(t)=a_\nu(t):=\left\{ \begin{array}{ll} -c & \hbox{if}\;\;
  t\in[0,\alpha) \cr
  b & \hbox{if}\;\; t\in[\alpha,1-\alpha] \cr 
  -\nu c & \hbox{if}\;\;
  t\in(1-\alpha,1] \end{array}\right.
\end{equation*}
with $\a\in(0,0.5)$, $b\geq 0$ and $c,\nu>0$. As a consequence, the problem is sublinear in some regions of the domain, while it is superlinear in others; for this reason it is referred to as \emph{superlinear indefinite}.

Observe that, when $\nu=1$, the weight $a(t)$ is symmetric in $\O$. For such choice, in \cite{LTZ}, it was proven that if $\l$ is sufficiently negative, there exists a value of $b$ such that problem \eqref{eq1} admits, for such choices of the parameters, an arbitrarily high number of solutions. Moreover, considering $b$ as main bifurcation parameter, the structure of the bifurcation diagrams was determined, showing that they consist of a \emph{primary branch} which possesses a certain number of turning points and contains symmetric solutions, and a number of \emph{secondary loops} bifurcating from it and containing asymmetric solutions. The analysis is based on the following observation: consider the set of all the positive solutions of the problem in the external intervals $[0,\a)$ and $(1-\a,1]$
\begin{equation}
\label{eq2}
  \!\!\text{(a) }\left\{ \begin{array}{l}
  -u''=\l u-c u^p \quad \hbox{in}\;\; (0,\a)\cr
  u(0)=M\end{array}\right. \quad
  \!\!\text{(b) }\left\{ \begin{array}{l}
  -u''=\l u-\nu c u^p \quad \hbox{in}\;\; (1-\a,1)\cr
  u(1)=M\end{array}\right.
\end{equation}
which will be denoted by $\mc{S}_0$ and $\mc{S}_{1,\nu}$ respectively. Then, if we define the following sets of points in the phase plane
\begin{equation}
\label{eq3}
\G_0:=\{(u(\a),u'(\a)): u\in\mc{S}_0\}, \qquad \G_{1,\nu}:=\{(u(1-\a),u'(1-\a)): u\in\mc{S}_{1,\nu}\},
\end{equation}
there is a $1-1$ correspondence between positive solutions of \eqref{eq1} and positive solutions of
\begin{equation}
\label{eq4}
  -u''=\l u+bu^p \quad \hbox{in}\;\; (\a,1-\a)
\end{equation}
which satisfy $(u(a),u'(\a))\in \G_0$ and $(u(1-\a),u'(1-\a))\in\G_{1,\nu}$. As a consequence, the study is reduced to consider in the phase plane associated to
\eqref{eq4} all the possible Poincar\'e maps that allow to connect $\G_0$ to $\G_{1,\nu}$ in a range of $1-2\a$, which is the amplitude of the central interval $(\a,1-\a)$ (observe that \eqref{eq4} is invariant by translations of the independent variable).

When $\nu\neq 1$, the weight $a(t)$ is no longer symmetric and the global structure of the bifurcation diagrams in $b$ of the solutions of \eqref{eq1} was determined in \cite{LT}, showing that they are not connected but exhibit an arbitrarily high number of bounded components, referred to as \emph{isolas}.
As $\nu$ perturbs from $1$ it is natural to expect that the secondary bifurcations occurring on the primary branch for $\nu=1$ transform into \emph{imperfect bifurcations}. Such phenomenon has not been considered in \cite{LT}, though it has been already detected numerically for problem \eqref{eq1} in \cite{LMT} and studied in general situations in \cite{LSW}. However, the sufficient conditions for imperfect bifurcations provided in \cite{LSW}, which are in the spirit of the functional-analytical transversality condition of the celebrated paper of M. G. Crandall and P. H. Rabinowitz \cite{CR}, cannot be tested in our situation, since the bifurcation points in \cite{LTZ} are constructed with the topological techniques described above. For this reason, we use here the same technique to prove the occurrence of imperfect bifurcations.

The work is distributed like follows: in Section \ref{section2} we recall some results of \cite{LTZ, LT} related to the geometry of the phase plane of the sublinear and superlinear problems; in Section \ref{section3} we provide the elements of \cite{LTZ} that allow to construct the bifurcation points in the symmetric case and in Section \ref{section4} we study the imperfect bifurcations in the asymmetric one.

\setcounter{equation}{0}
\section{The geometry of the phase plane}
\label{section2}
We start by recalling from \cite[Section 2]{LT} the geometry of the flow of the sublinear problems \eqref{eq2}, describing the set of points in the phase plane $(u,u')$ reached at time $\a$ and $1-\a$ respectively.

\begin{figure}[ht]
\begin{center}
\begin{tabular}{ccc}
\includegraphics[scale=0.50]{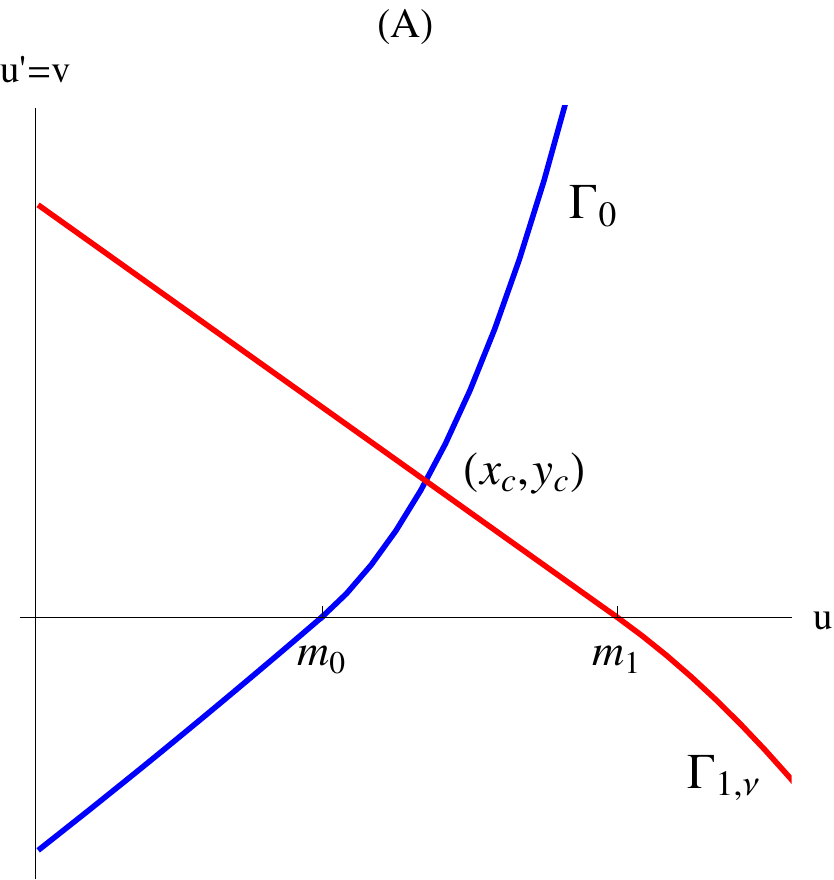} & 
\includegraphics[scale=0.50]{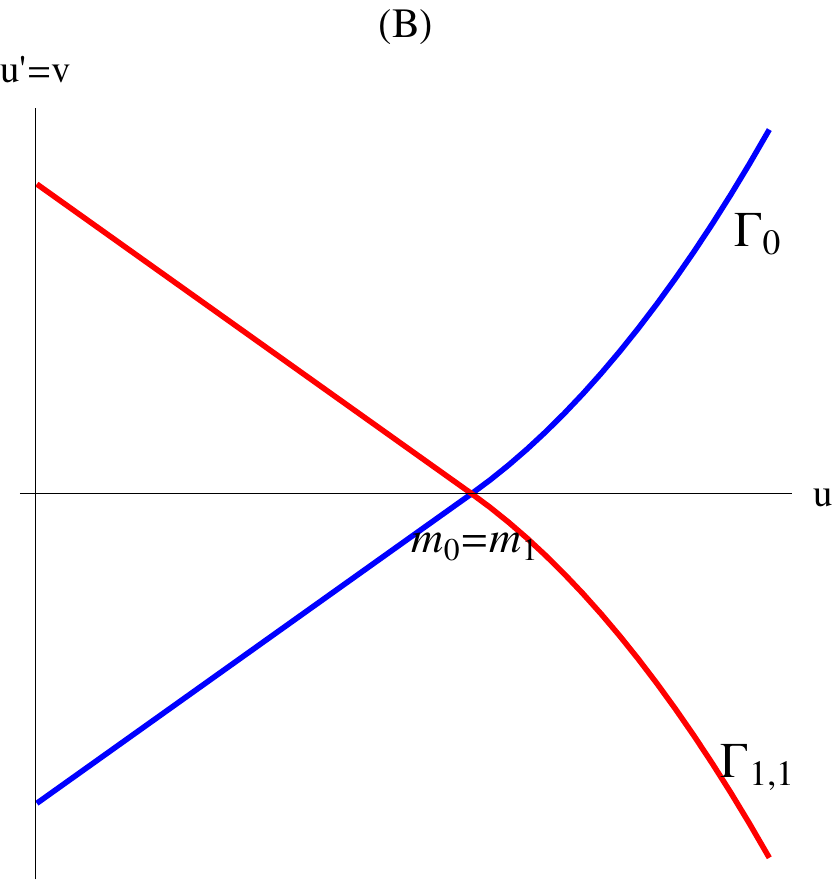} &
\includegraphics[scale=0.50]{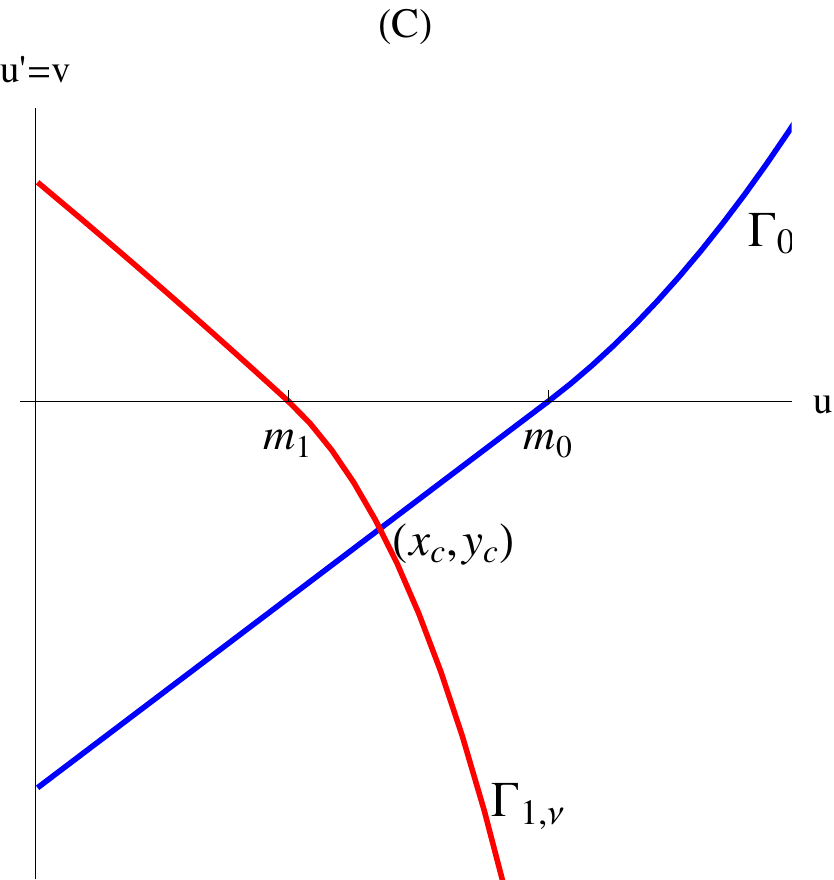}
\end{tabular}\\
\caption{Relative positions of the curves $\G_0$ and $\G_{1,\nu}$ for $\nu<1$ (A), $\nu=1$ (B) and $\nu>1$ (C).} \label{Fig1}
\end{center}
\end{figure}

\begin{proposition}
\label{pr21}
For every $\l\leq 0$, $M\in(0,\infty]$, $c>0$ and $\a\in(0,0.5)$ there exists a unique
increasing function $y_0^c\in\mathcal{C}^1([0,\infty),\R)$ such that $\lim_{x\ua +\infty}y_0^c(x)=+\infty$ and
\begin{itemize}
\item[\rm (i)] $u'(\a)=y_0^c(u(\a))$ for any nonnegative solution $u$ of \eqref{eq2}(a);
\item[\rm (ii)] there is a unique $m_0>0$ such that $y_0^c(m_0)=0$;
\item[\rm (iii)] if $\tilde{c}>c$ and $x\in[0,+\infty)$, then $y_0^{\tilde{c}}(x)>y_0^{c}(x)$.
\end{itemize}
\end{proposition}
Performing the change of variable $\tilde t=1-t$, which entails a change of sign of the first derivatives of the solutions, we obtain immediately a counterpart of Proposition \ref{pr21} for Problem \eqref{eq2}(b). As a consequence, maintaining the notation of \eqref{eq3}, the sets $\G_0$ and $\G_{1,\nu}$ are differentiable curves whose relative position, thanks to Proposition \ref{pr21}(iii), is like shown in Figure \ref{Fig1}. In particular, they intersect in a point $(x_c,y_c)$ such that $\sign(y_c)=\sign(1-\nu)$. 

Passing to the phase plane associated to the flow in the central interval, observe that equation \eqref{eq4} possesses two nonnegative equilibria: $(0,0)$ and $(\O,0)$, where
\begin{equation}
\label{eq5} \O=\left(-\l/b\right)^{\frac{1}{p-1}}.
\end{equation}
Moreover, the first order system associated to \eqref{eq4} admits the
first integral
\begin{equation}
\label{eq6}
  E(u,v):=v^2+\lambda u^2+\frac{2b}{p+1}u^{p+1}
\end{equation}
whose associated potential energy is $\v(u):=\lambda u^2+\frac{2b}{p+1}u^{p+1}$. As $\v''(0)=2\l<0$ and $\v''(\Omega)=2\l(1-p)>0$, $(0,0)$ is a
saddle point and $(\O,0)$ is a center. Moreover, as $\v(u)$ has the shape of a  potential  well, 
there is an homoclinic connection of $(0,0)$ which will be denoted by $\g_h=\g_h(b)$ and
surrounds $(\O,0)$ and every  closed orbit around it. Its equation, thanks to \eqref{eq6}, is
\begin{equation}
\label{eq7}
v(u)=\pm\sqrt{-\l u^2 - \frac{2b}{p+1}u^{p+1}}, \qquad 0\leq u\leq\left(\frac{-\l(p+1)}{2b}\right)^{\frac{1}{p-1}}
\end{equation}
and, therefore, it approaches the lines $v=\pm\sqrt{-\l} u$ as $b\da 0$ and shrinks monotonically to $(0,0)$ as $b\ua+\infty$.

Observe that, if we set
\begin{equation}
\label{eq8}
b^*=b^*(\l)=-\l/m_0^{p-1}(\l),
\end{equation}
where $m_0(\l)$ is the one of Proposition \ref{pr21}(ii), then, in view of \eqref{eq5}, $\O=m_0$ if $b=b^*$, while $\O>m_0$ if $b<b^*$ and $\O<m_0$ if $b>b^*$. In this latter case, if we superimpose the curves $\G_0$ and $\G_{1,1}$ to the flow of \eqref{eq4}, we have that small closed orbits surrounding $(\O,0)$ do not touch $\G_0$ and $\G_{1,1}$, while orbits far away from it do. As a consequence, there will be an integral line of the flow which intersects the curves for the first time. We will denote it by $\g_t=\g_t(b)$ and assume that it the unique orbit that is tangent to the curves $\G_0$ and $\G_{1,1}$, while larger integral lines are secant to them; the abscissa of the tangency point will be denoted by $x_{t,0}=x_{t,0}(b)$. As shown in \cite[Appendix]{LT}, this is the case at least if $M<+\infty$ and $c\sim 0$. We make this assumption in order to simplify the discussion of the subsequent sections; however, should it not be satisfied, this would only turn into more complex bifurcation diagrams (with more bifurcation or turning points or components) and our analysis can be easily adapted.

With this assumption, we have in particular that there exists a value of $b$, denoted by $b_h$, for which the homoclinic $\g_h$ introduced in \eqref{eq7} is tangent to $\G_0$ and $\G_{1,1}$ and does not intersect them for $b>b_h$ (see Figure \ref{Fig2}(A)).

\setcounter{equation}{0}
\section{Construction of the bifurcation point in the symmetric case}
\label{section3}
We start this section recalling the main results related to multiplicity of solutions and the structure of the bifurcation diagrams in the symmetric case $\nu=1$. They have been obtained in \cite{LTZ} and can be summarized as follows.

\begin{theorem}
\label{th22}
(i) Fix $\nu=1$ and set, for $j\in\N\setminus\{0\}$, $\l_j:=\frac{1}{1-p}\left(\frac{2\pi j}{1-2\a}\right)^2<0$.
Then, for every $c>0$, $p>1$, $M\in (0,\infty]$ and $\l\in (\l_{j+1},\l_{j})$, if we take $b=b^*$ as in \eqref{eq8}, Problem \eqref{eq1} admits, at least, $4j$ solutions.
\par
(ii) With the same choice of the parameters, the minimal global bifurcation diagram in $b$ of the positive solutions of Problem \eqref{eq1} consists of a \emph{primary curve}, filled in by symmetric solutions, plus $j$ \emph{secondary loops}, filled in by asymmetric solutions and emerging each one at two bifurcations points from the primary curve which will be denoted by $b_b^{i,-}$ and $b_b^{i,+}$ in such a way that
\begin{equation}
\label{eq9}
b_b^{i,-}<b_b^{i+1,-}<b^*<b_b^{i+1,+}<b_b^{i,+} \qquad \text{for $1\leq i\leq j-1$.}
\end{equation}
\end{theorem}

In the rest of this section we will sketch some elements of the proof of Theorem \ref{th22}, up to arrive to the construction of the bifurcation points, that will be fundamental in the next section. Using the notation introduced in \eqref{eq9}, we will detail the case $b=b_b^{1,+}$. The other ones can be treated analogously, with minor changes.

As already commented in Section \ref{section1}, the key ingredient in our method is the superposition of the curves $\G_0$ and $\G_{1,1}$ to the flow of equation \eqref{eq4}, in order to study how to connect the former curve to the latter along the flow in a range of $1-2\a$. For this reason we introduce the Poincar\'e maps $\tau_j(x,b)$, $j\geq 1$, where, for every $x\neq x_{t,0}$ in the appropriate domain, $\tau_j(x,b)$ denotes the range necessary to reach the curve $\G_{1,1}$ for the $j-$th time, starting from $(x,y^c_0(x))\in\G_0$ and moving along the flow of \eqref{eq4}, remaining in the region $u>0$. We do not define the Poincar\'e maps for $x= x_{t,0}$ at this stage, since in that case the orbit of the flow is tangent to $\G_{1,1}$ and the intersection points must be counted with multiplicity; we will cover this case later.

As a consequence of the discussion carried out in Section \ref{section1}, we have that every
\begin{equation*}
x\in U(b):=\bigcup_{j\geq 1}\tau_j(\cdot,b)^{-1}(1-2\a)
\end{equation*}
provides us with a different solution of Problem \eqref{eq1} with $\nu=1$.

\begin{figure}[ht]
\begin{center}
\begin{tabular}{cc}
\includegraphics[scale=0.50]{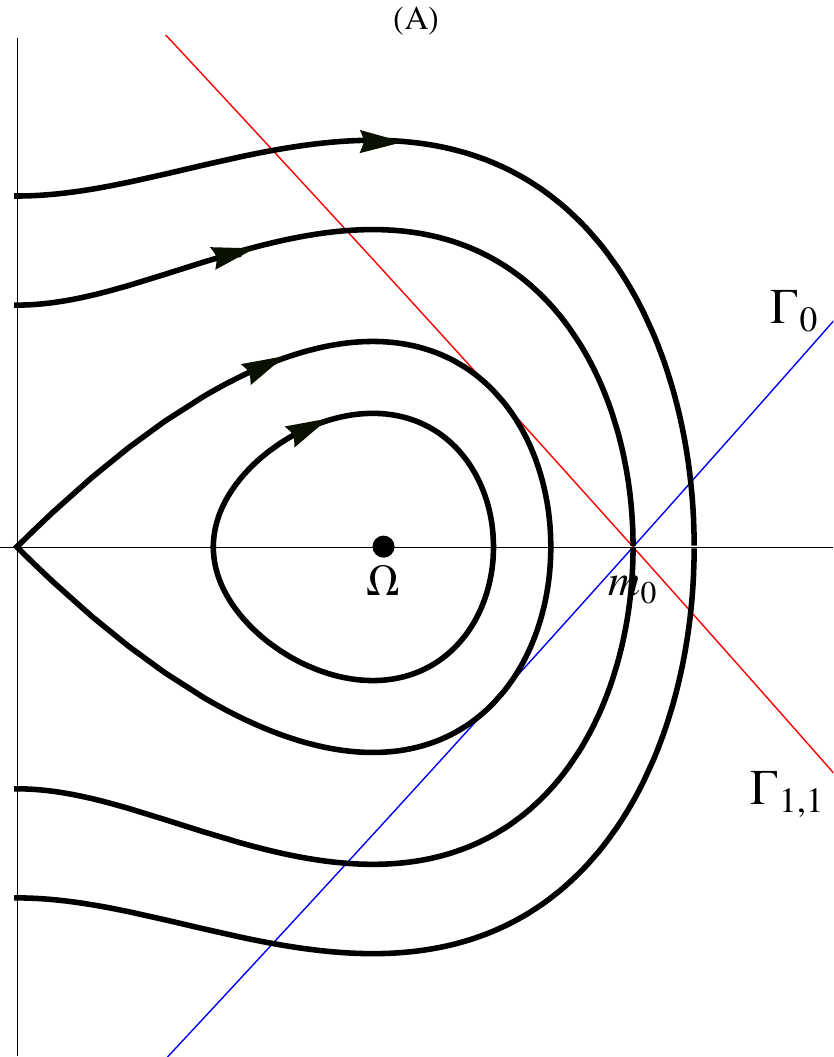} & \hspace{0.5cm}
\includegraphics[scale=0.6]{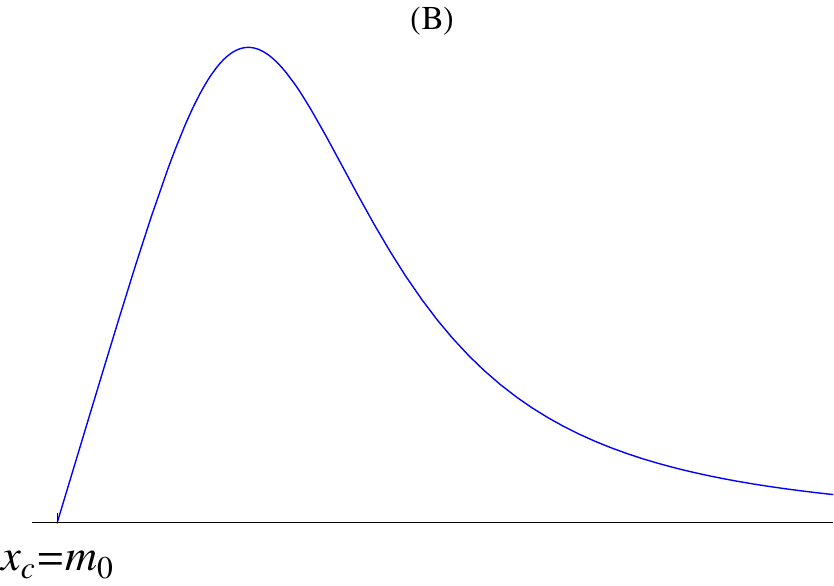} 
\end{tabular}\\
\caption{Phase plane for $b=b_h$ (A) and $\tau_1(\cdot,b)$, the only Poincar\'e map defined for $b\geq b_h$ (B).} \label{Fig2}
\end{center}
\end{figure}

Our first claim in the construction of the bifurcation points is that, at such points, $b$ must be smaller than $b_h$. Indeed, if $b\geq b_h$, the unique Poincar\'e map that can be defined is $\tau_1(x,b)$ for $x>m_0$, which is a regular function of $x$ and $b$ (see Figure \ref{Fig2}). As a consequence, no bifurcation point can occur.

\begin{figure}[ht]
\begin{center}
\begin{tabular}{ccc}
\includegraphics[scale=0.50]{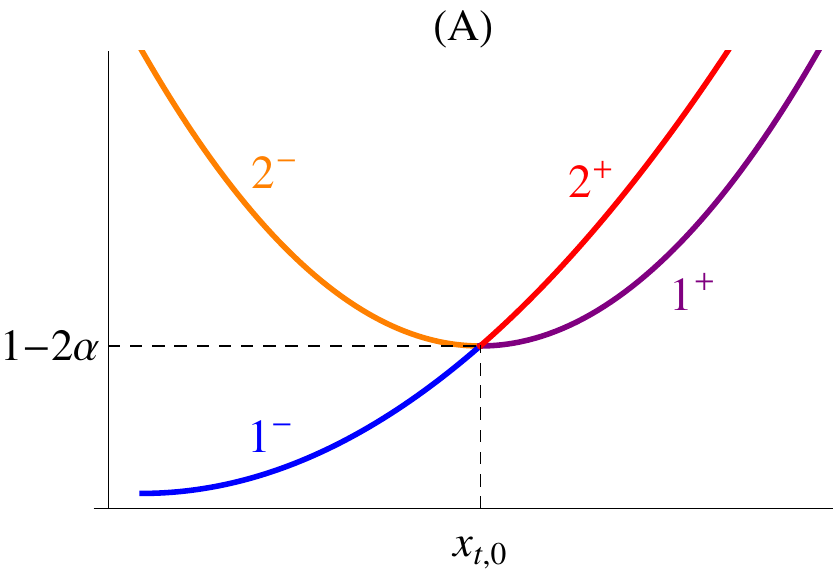} &  
\includegraphics[scale=0.50]{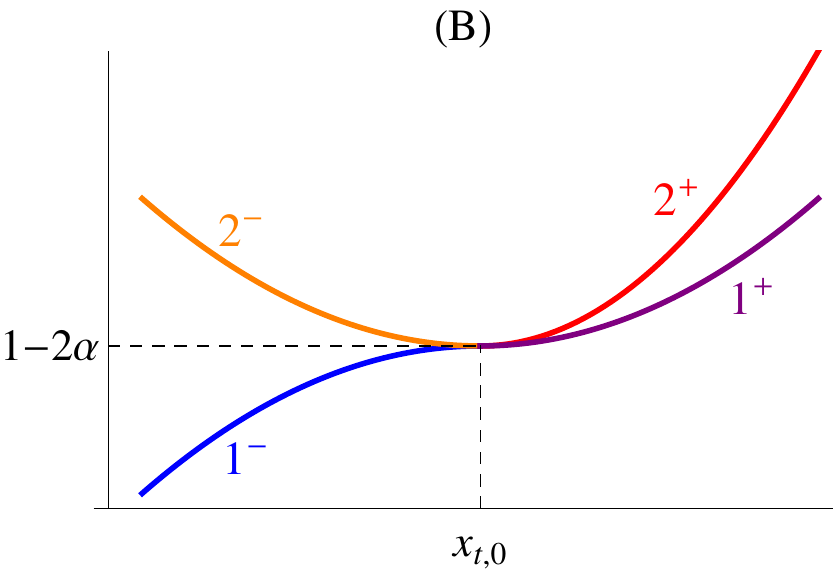} & 
\includegraphics[scale=0.50]{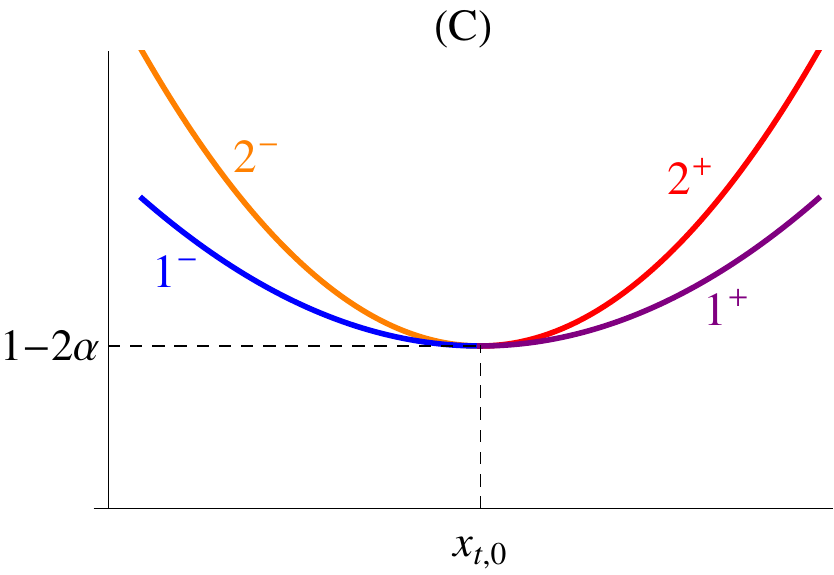}
\end{tabular}
\caption{Possible configurations of the Poincar\'e maps $\t_1(\cdot,b)$ and $\t_2(\cdot,b)$.} \label{Fig3}
\end{center}
\end{figure}

Therefore we will take $b<b_h$ and, from the discussion of Section \ref{section2}, we have that $\g_h$ intersects $\G_0$ in two points, whose abscissas will be denoted by $x_{h,0}^-$ and $x_{h,0}^+$. For the same reason, all the closed orbits between $\g_t$ and $\g_h$ intersect $\G_{1,1}$ in two points. We denote the orbit which passes through $(x,y_0^c(x))\in\G_0$ by $\g_x$; thanks to \eqref{eq6} its equation will be
\begin{equation}
\label{eq10}
v^2+\l u^2+\frac{2b}{p+1}u^{p+1}=(y_0^c(x))^2+\l x^2+\frac{2b}{p+1}x^{p+1}.
\end{equation}
Moreover, if we denote the abscissas of its intersection points with $\G_{1,1}$ by $x_1^-(x)$ and $x_1^+(x)$, in such a way that
\begin{equation}
\label{eq11}
x_1^-(x)<x_1^+(x),
\end{equation}
it is easy to see that they satisfy
\begin{equation}
\label{eq12}
\lim_{x\to x_{t,0}} x_1^-(x)= x_{t,0} =\lim_{x\to x_{t,0}} x_1^+(x).
\end{equation}
We also denote by $x_0(x)$ the abscissa of the point of intersection between $\G_0$ and $\g_x$ which is different from $(x,y_0^c(x))$. With these notations, thanks to the symmetry of the problem for $\nu=1$, we have that
\begin{equation}
\label{eq13}
  \left. \begin{array}{cc}
x_1^-(x)=x \quad \text{ and } \quad x_0(x)=x_1^+(x) & \text{ if $x<x_{t,0}$,} \\
x_1^+(x)=x \quad \text{ and } \quad x_0(x)=x_1^-(x) & \text{ if $x>x_{t,0}$.} 
\end{array}\right.
\end{equation}
In addition, $\g_x$ intersects the $u-$axis in two points, which will be denoted by $(x_m(x),0)$ and $(x_M(x),0)$ with $x_m(x)<x_M(x)$ and such that
\begin{equation}
\label{eq14}
\lim_{x\to x_{t,0}} x_m(x)= x_m(x_{t,0}), \qquad \lim_{x\to x_{t,0}} x_M(x)= x_M(x_{t,0}).
\end{equation}
If we set
\begin{equation*}
I(u,x,b):=\left((y_0^c(x))^2+\l \left(x^2-u^2\right)+\frac{2b}{p+1}\left(x^{p+1}-u^{p+1}\right)\right)^{-\frac{1}{2}},
\end{equation*}
we have from our definitions and \eqref{eq10} that, for $x\sim x_{t,0}$, $x\neq x_{t,0}$,
\begin{align*}
\tau_1(x,b)&= \int_{x_m(x)}^{x}I(u,x,b)\,du + \int_{x_m(x)}^{x_1^-(x)}I(u,x,b)\,du,
\\
\tau_2(x,b)&= \int_{x_m(x)}^{x}I(u,x,b)\,du + \int_{x_m(x)}^{x_1^+(x)}I(u,x,b)\,du.
\end{align*}
From \eqref{eq11} it is immediate that
\begin{equation}
\label{eq15}
\tau_1(x,b)<\tau_2(x,b) \quad \text{ for $x\sim x_{t,0}$, $x\neq x_{t,0}$}
\end{equation}
and, in view of \eqref{eq13},
\begin{align}
\label{eq16}
\t_1(x,b)&:=2 \int_{x_m(x)}^{x}\!\!I(u,x,b)\,du=\begin{cases}
\tau_1(x,b) & \text{ if $x<x_{t,0}$} \\
\tau_2(x,b) & \text{ if $x>x_{t,0}$} 
\end{cases},
\\
\label{eq17}
\t_2(x,b)&:= \int_{x_m(x)}^{x}\!\!I(u,x,b)\,du + \int_{x_m(x)}^{x_0(x)}\!I(u,x,b)\,du=\begin{cases}
\tau_2(x,b) & \text{ if $x<x_{t,0}$} \\
\tau_1(x,b) & \text{ if $x>x_{t,0}$} 
\end{cases}.
\end{align}
Thanks to \eqref{eq12} and \eqref{eq14}, it is possible to extend these maps continuously, by setting
\begin{equation*}
\t_1(x_{t,0},b)=\t_2(x_{t,0},b):=\lim_{x\to x_{t,0}}\tau_1(x,b)=\lim_{x\to x_{t,0}}\tau_2(x,b)
\end{equation*}
and, since $x_m(x)$ and $x_0(x)$ are analytic functions of $x$, it follows from \eqref{eq16} and \eqref{eq17} that the maps $\t_j(\cdot,b)$, $j\in\{1,2\}$ also are analytic. As $x\da x_{h,0}^-$ the orbit $\g_x$ approaches the equilibrium $(0,0)$ and, by continuous dependence, this entails
\begin{equation*}
\lim_{x\da x_{h,0}^-}\t_1(x,b)=+\infty=\lim_{x\da x_{h,0}^-}\t_2(x,b).
\end{equation*}
In particular we obtain that these maps are non-constant and that there exists a lowest-order derivative (in $x$) of them which does not vanish. We can be more precise for the case of $\t_2$: indeed, from \eqref{eq17} and \eqref{eq13}, we also have that
\begin{equation}
\label{eq18}
\t_2(x,b)=\t_2(x_0(x),b)
\end{equation}
and, since $x_0(\cdot)$ is decreasing due to the geometry of the phase plane, by successively differentiating \eqref{eq18} and using \eqref{eq12} and \eqref{eq14}, it is possible to show that the lowest-order derivative of $\t_2(\cdot,b)$ which does not vanish, when evaluated in $x_{t,0}$, is of even order. As a consequence, recalling \eqref{eq15}, we obtain that $\p\t_1(x_{t,0}(b),b)/\p x\geq 0$. Summarizing, we have that the four branches of Poincar\'e maps ($1^-$) $\tau_1(x,b)$ for $x<x_{t,0}$, ($2^-$) $\tau_2(x,b)$ for $x<x_{t,0}$, ($1^+$) $\tau_1(x,b)$ for $x>x_{t,0}$ and ($2^+$) $\tau_2(x,b)$ for $x>x_{t,0}$ meet at $x=x_{t,0}$ in one of the possible configurations represented in Figure \ref{Fig3} or in analogous ones, taking into account all the possibilities for their concavities.

\begin{figure}[ht]
\begin{center}
\begin{tabular}{ccc}
\includegraphics[scale=0.50]{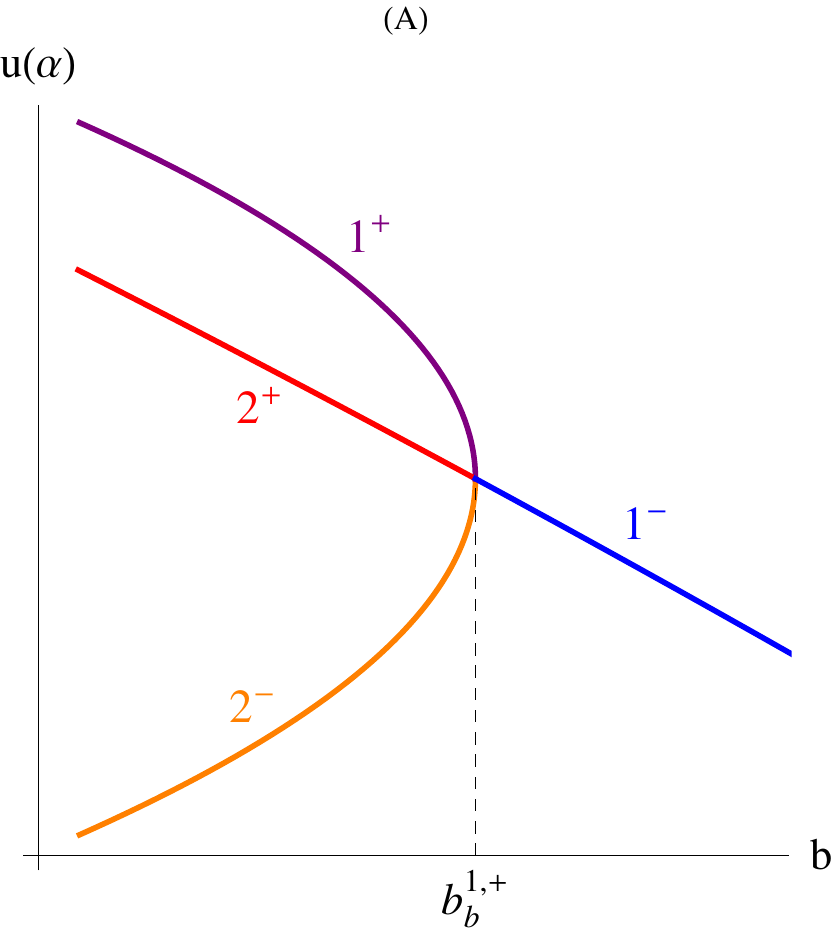} & 
\includegraphics[scale=0.50]{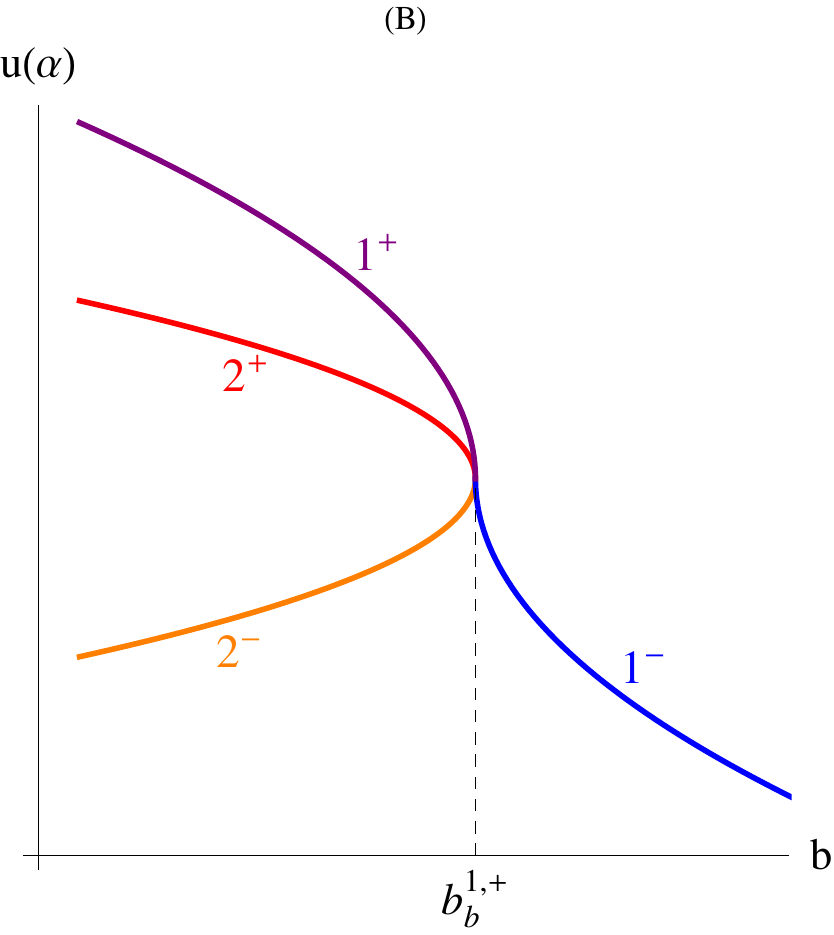} & 
\includegraphics[scale=0.50]{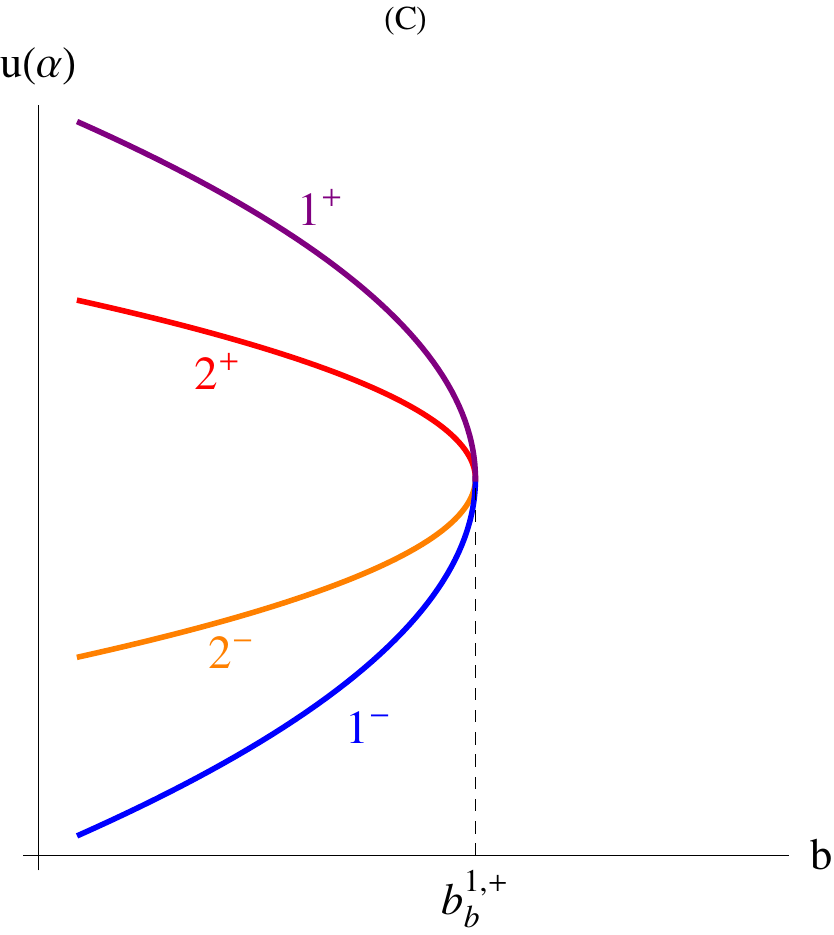}
\end{tabular}\\
\caption{Possible structure of the bifurcation points for $\nu=1$.} \label{Fig4}
\end{center}
\end{figure}

By linearizing \eqref{eq4} at $\O=m_0$ for $b=b^*$, we obtain that the period of small oscillations around $(\O,0)$ is given by
\begin{equation*}
\tau:=\frac{2\pi}{\sqrt{-\l(p-1)}}
\end{equation*}
and therefore, with the notation of Theorem \ref{th22}, if $\l<\l_1$, it follows that
\begin{equation}
\label{eq19}
\lim_{x\to m_0}\tau_2(x,b^*)<\tau<1-2\a.
\end{equation}
Intuitively, since $\g_{t}(b)$ shrinks towards $(m_0,0)$ as $b\da b^*$, we have that $x_{t,0}(b)\to m_0$ as $b\da b^*$ and
\begin{equation*}
\lim_{b\da b^*}\tau_2(x_{t,0}(b),b)=\lim_{x\to m_0}\tau_2(x,b^*)
\end{equation*}
(see \cite[Theorem 4.2]{LTZ} for the rigorous proof, which is quite delicate, since it requires a singular perturbation argument). This, together with \eqref{eq19}, entails
\begin{equation}
\label{eq20}
\tau_2(x_{t,0}(b),b)<1-2\a \quad \text{ for $b>b^*$, $b\sim b^*$.}
\end{equation}
On the other hand, as $b\ua b_h$, the domains of the Poincar\'e maps related to the closed orbits shrink towards $x_{t,0}(b_h)$ and such Poincar\'e maps blow up to $+\infty$, as a consequence of the continuous dependence theorem, since, moving along the flow of \eqref{eq4}, we have to pass in a vicinity of the equilibrium $(0,0)$; therefore
\begin{equation*}
\tau_2(x_{t,0}(b),b)>1-2\a \quad \text{ for $b<b_h$, $b\sim b_h$}.
\end{equation*}
Combining this with \eqref{eq20} and again the continuous dependence of the time maps with respect to $b$, we obtain that
\begin{multline}
\label{eq21}
B^{1,+}:=\{\ov{b}>b^*: \tau_2(x_{t,0}(b),b)<1-2\a \text{ for $b<\ov{b}$, $b\sim\ov{b}$, } \tau_2(x_{t,0}(\ov{b}),\ov{b})=1-2\a \\ 
\text{ and } \tau_2(x_{t,0}(b),b)>1-2\a \text{ for $b>\ov{b}$, $b\sim\ov{b}$} \}\neq\emptyset.
\end{multline}
As a consequence, we have that
\begin{equation}
\label{eq22}
b_b^{1,+}=\max B^{1,+}.
\end{equation}
Since $\t_2(\cdot,b)$ has a minumum (or a maximum) at $x_{t,0}$, the branches $(2^-)$ and $(1^+)$ form a subcritical (or supercritical) turning point at $b=b_b^{1,+}$, while different situations can occur according to the behavior of $\t_1(x,b)$ at $x=x_{t,0}$ and $b=b_{b}^{1,+}$: either it is transversal to $\t_2$ (see Figure \ref{Fig3}(A)) or it is strictly increasing in a neighborhood of $x_{t,0}$ without being transversal to $\t_2$ (see Figure \ref{Fig3}(B)), or it has the same monotony as $\t_2$ (see Figure \ref{Fig3}(C)). According to each of these situations, we obtain a different configuration of the bifurcation point: respectively transcritical non degenerate pitchfork bifurcation, transcritical degenerate pitchfork bifurcation and double pitchfork bifurcation (see Figure \ref{Fig4}(A),(B) and (C) respectively, where, as in all the subsequent bifurcation diagram, we represent the value of $b$ on the abscissas versus the value of the solutions at $t=\a$; the numbers refer to the branch of Poincar\'e map of Figure \ref{Fig3} that generates the corresponding branch of the bifurcation point).

\setcounter{equation}{0}
\section{Imperfect bifurcation}
\label{section4}
Now we perturb the symmetric case by taking $\nu\neq 1$, $\nu\sim 1$ and study how the bifurcation diagrams change in a neighborhood of a bifurcation point. As our analysis is strongly based on the previous sections, we will often refer to the constructions introduced there and maintain the same notation. 

\subsection{The case $\nu>1$} 
\begin{figure}[ht]
\begin{center}
\includegraphics[scale=0.68]{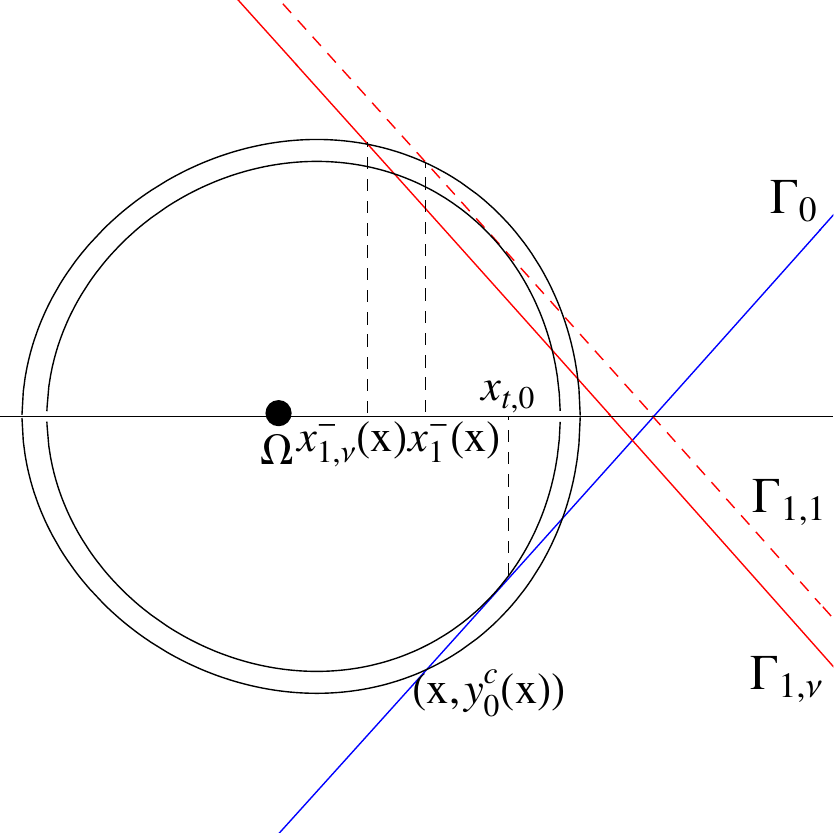} 
\caption{Geometry of the phase plane for $\nu>1$ in a neighborhood of $\g_t$.} \label{Fig5}
\end{center}
\end{figure}

In this case, thanks to the analysis of Section \ref{section2} (see Figure \ref{Fig1}(C)), the geometry of the phase plane, in a neighborhood of the tangent orbit to $\G_0$, looks like represented in Figure \ref{Fig5}.

In particular, any orbit $\g_x$ with $x\sim x_{t,0}$ intersects $\G_{1,\nu}$ in two points whose abscissas will be denoted by $x_{1,\nu}^-(x)$ and $x_{1,\nu}^+(x)$ in such a way that
\begin{equation}
\label{eq23}
x_{1,\nu}^-(x)<x_{1}^-(x)\leq x_{t,0}\leq x_{1}^+(x)<x_{1,\nu}^+(x).
\end{equation}
The main differences with respect to the case of Section \ref{section3} are that such orbits $\g_x$ are now separated away from the orbit $\g_{t}^{1,\nu}$, the tangent one to $\G_{1,\nu}$ and that $\g_t$ itself is now secant to $\G_{1,\nu}$.

We introduce the perturbed Poincar\'e maps, $\tau_{j,\nu}(x,b)$, $j\geq 1$, for every $x\sim x_{t,0}$ in the appropriate domain, as the time necessary to reach the curve $\G_{1,\nu}$ for the $j-$th time, starting from $(x,y^c_0(x))\in\G_0$ and moving along the flow of \eqref{eq4} remaining in the region $u>0$. Similarly as before, from this definition, we have that
\begin{align*}
\tau_{1,\nu}(x,b)&= \int_{x_m(x)}^{x}I(u,x,b)\,du + \int_{x_m(x)}^{x_{1,\nu}^-(x)}I(u,x,b)\,du,
\\
\tau_{2,\nu}(x,b)&= \int_{x_m(x)}^{x}I(u,x,b)\,du + \int_{x_m(x)}^{x_{1,\nu}^+(x)}I(u,x,b)\,du.
\end{align*}
Thanks to these relations and \eqref{eq23}, it is immediate that for every $x\sim x_{t,0}$, $b\sim b_{b}^{1,+}$,
\begin{equation}
\label{eq24}
\tau_{1,\nu}(x,b)<\tau_{1}(x,b)\leq\tau_{2}(x,b)<\tau_{2,\nu}(x,b).
\end{equation}
Moreover now, since the orbits $\g_x$ for $x\sim x_{t,0}$ are separated from $\g_{t}^{1,\nu}$, there is no multiple intersection with $\G_{1,\nu}$ and, therefore, the two branches $\tau_{j,\nu}(x,b)$, $j\in\{1,2\}$, which are differentiable for every $x\sim x_{t,0}$, $b\sim b_{b}^{1,+}$ by the differentiable dependence theorem, are also separated, contrary to what happens in Section \ref{section3}.

With these elements we are now able to present the main result of this paper. In it we make an additional assumption which essentially says that the bifurcation point is non-degenerate; we will do some remarks on the general case afterwards.

\begin{theorem}
Assume that
\begin{equation}
\label{eq25}
\frac{\p^2 \t_j(x_{t,0},b_{b}^{1,+})}{\p x^2}\neq 0, \quad \text{ $j\in\{1,2\}$}.
\end{equation}
Then there exists $\ov{\nu}>1$, $\ov{\nu}\sim 1$ and a neighborhood of $b_b^{1,+}$ such that, if $1<\nu<\ov{\nu}$, the bifurcation point occurring for $\nu=1$ and $b=b_b^{1,+}$ gives rise to an imperfect bifurcation of one of the types represented in Figure \ref{Fig6}.

\begin{figure}[h]
\begin{center}
\begin{tabular}{cc} 
\includegraphics[scale=0.50]{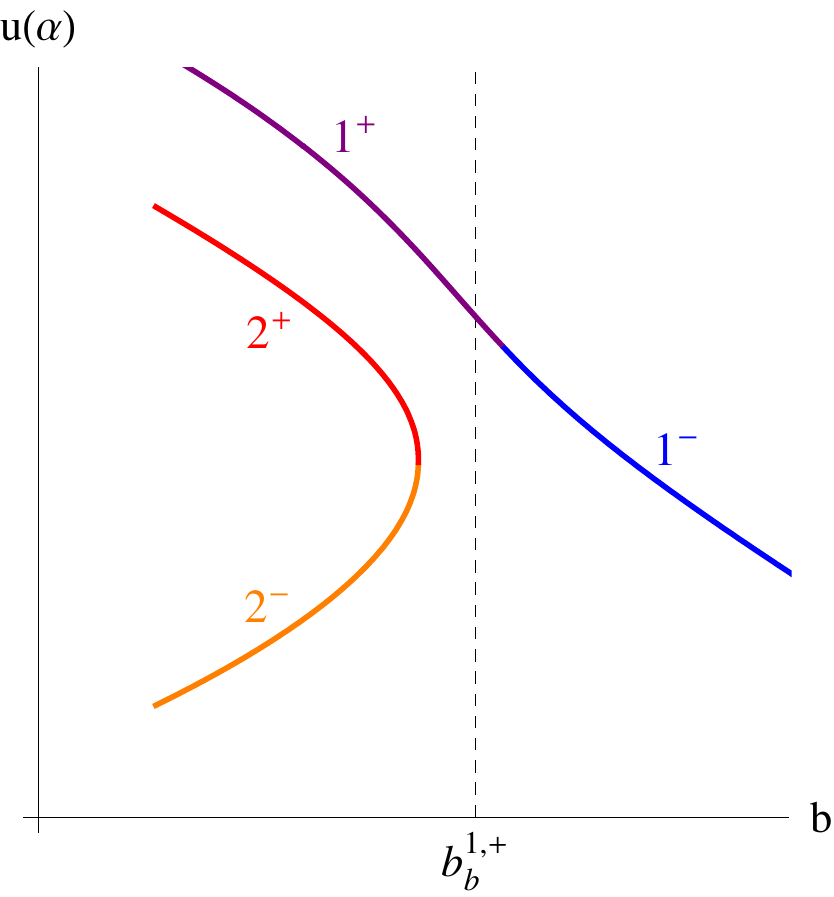} & \hspace{0.5cm}
\includegraphics[scale=0.50]{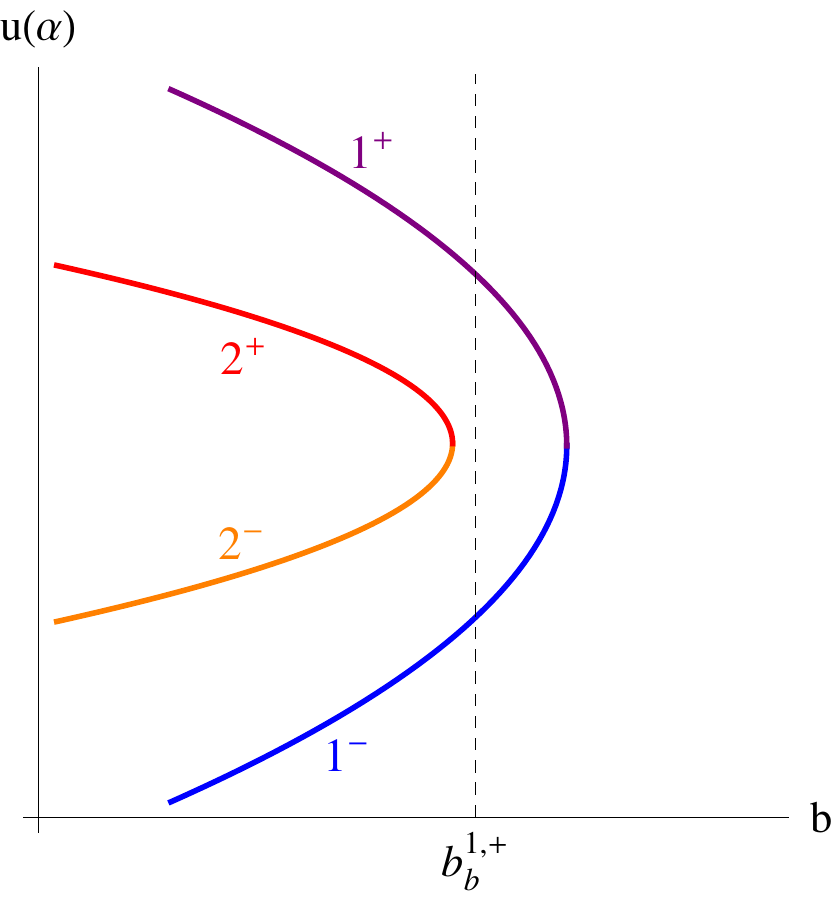}
\end{tabular}\\
\caption{Possible structure of the imperfect bifurcation under condition \eqref{eq25}.} \label{Fig6}
\end{center}
\end{figure}

\begin{proof}
We start with the case corresponding to Figure \ref{Fig3}(A). Recalling the definition of $b_b^{1,+}$ (see \eqref{eq21} and \eqref{eq22}), there exists $\d> 0$, $\d\sim 0$ such that, for $b\in(b_b^{1,+}-2\d,b_b^{1,+}-\d)$ and $x\sim x_{t,0}$, $\tau_2(x,b)<1-2\a$. By continuous dependence, there exists $\ov{\nu}>1$, $\ov{\nu}\sim 1$ such that, if $1<\nu<\ov{\nu}$, $\tau_{2,\nu}(x,b)<1-2\a$ for the same $x$'s and $b$'s. However, \eqref{eq24} shows that $\tau_{2,\nu}(x,b_b^{1,+})>1-2\a$ for $x\sim x_{t,0}$. Moreover, by the differentiable dependence theorem, the derivative with respect to $x$ of $\tau_{j,\nu}(x,b)$, $j\in\{1,2\}$, is as close as we wish, if $\nu\sim 1$, to the one of $\tau_{j,\nu}(x,b)$, for $x\sim x_{t,0}$, $x\neq x_{t,0}$. In particular, reducing $\ov{\nu}$ and $\d$ if necessary, $\tau_{2,\nu}(x,b)$ is decreasing in $\mc{U}^-:=(x_{t,0}-2\d,x_{t,0}-\d)$ and increasing in $\mc{U}^+:=(x_{t,0}+\d,x_{t,0}+2\d)$. Therefore, for $\nu\sim 1$, $\tau_{2,\nu}(x,b)$ has a minimum in a neighborhood of $x_{t,0}$, which is unique thanks to \eqref{eq25}. Similarly, $\tau_{1,\nu}(x,b)$ is increasing in $\mc{U}^-\cup\mc{U}^+$ and by \eqref{eq25} it is increasing in a whole neighborhood of $x_{t,0}$. Moreover, by continuous dependence and \eqref{eq24}, $\tau_{1,\nu}(x,b)<1-2\a$ for $b\sim b_{b}^{1,+}$ and $x\in\mc{U}^-$, while $\tau_{1,\nu}(x,b)>1-2\a$ for $x\in\mc{U}^+$. Thus, for each $b\sim b_{b}^{1,+}$, there exists a unique value of $x(b)$ such that $\tau_{1,\nu}(x(b),b)=1-2\a$. Summarizing, the perturbed Poincar\'e maps look like represented in the left graph of Figure \ref{Fig7}.

\begin{figure}[h]
\begin{center}
\begin{tabular}{cc}
\includegraphics[scale=0.55]{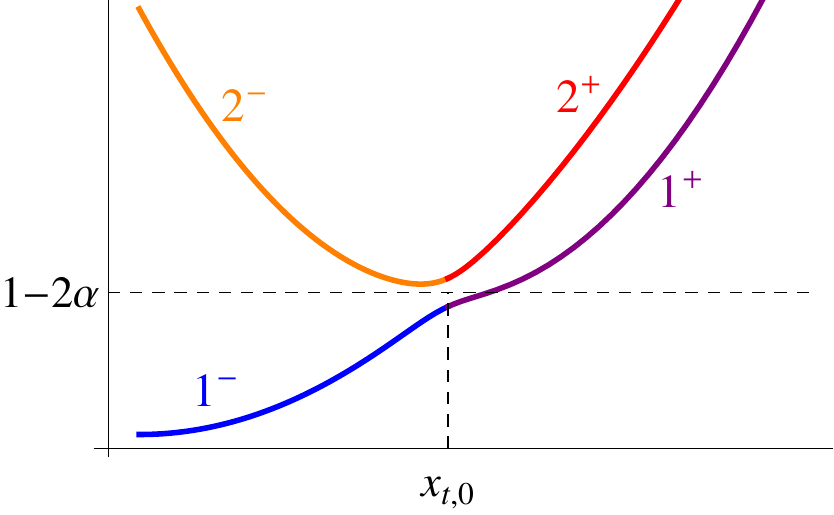} & \hspace{0.5cm}
\includegraphics[scale=0.55]{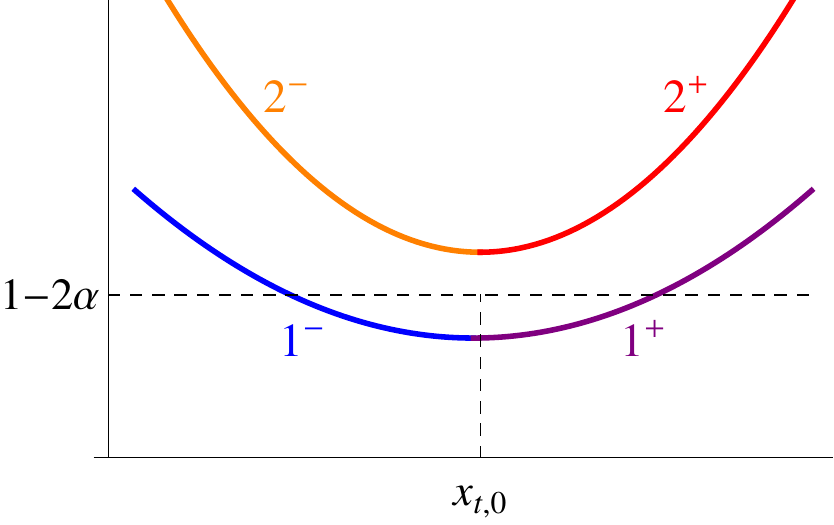}
\end{tabular}\\
\caption{Perturbed Poincar\'e maps for $b=b_b^{1,+}$, corresponding to the case of Figure \ref{Fig3}(A) (left) and \ref{Fig3}(C) (right).} \label{Fig7}
\end{center}
\end{figure}

From the previous argument, we have that, for $1<\nu<\ov{\nu}$, $\tau_{2,\nu}(\cdot,b)$ produces a subcritical turning point for a certain value of $b$ smaller than $b_b^{1,+}$, precisely the one for which its minimum crosses the value $1-2\a$, while $\tau_{1,\nu}(\cdot,b)$ produces a regular branch, as shown in the left diagram of Figure \ref{Fig6}.

Reasoning similarly, in the case corresponding to Figure \ref{Fig3}(C) $\tau_{1,\nu}(\cdot,b)$ also has a minimum in a neighborhood of $x_{t,0}$ which crosses the value $1-2\a$ for a certain $b>b_b^{1,+}$, since, as in section \ref{section3}, it blows to $+\infty$ as $b\ua b_h$. Therefore, it produces another turning point for $b>b_b^{1,+}$. This situation corresponds to the right pictures in Figures \ref{Fig7} and \ref{Fig6}.
\end{proof}
\end{theorem}

When \eqref{eq25} does not hold, which is the case, for example, of Figure \ref{Fig3}(B), a deeper analysis than in the previous proof is required, since the zero of the derivative of $\t_2$ is not simple and it could perturb, as $\nu\neq 1$, in more than one zero of $\tau_{j,\nu}'(\cdot,b)$, $j\in\{1,2\}$. In the bifurcation diagrams this would lead to the presence of more turning points on each of the two branches of the imperfect bifurcation. This analysis goes outside the scope of this work, therefore we stop here.

\subsection{The case $\nu<1$} If we perform the change of variable $\tilde{u}(t):=u(1-t)$, then $u$ is a solution of Problem \eqref{eq1} if and only if $\tilde{u}$ solves
\begin{equation}
\label{eq26}
  \left\{ \begin{array}{l}
  -\tilde{u}''=\l \tilde{u}+\tilde{a}(t)\tilde{u}^p \quad \hbox{in}\;\; (0,1)\cr
  \tilde{u}(0)=\tilde{u}(1)=M,\end{array}\right.
\end{equation}
where we have set $\tilde{a}(t):=a(1-t)$. Now observe that, if $t\in(0,\a)$,
\begin{equation*}
\tilde{a}(t)=a(1-t)=-\nu c>-c=a(t)=\tilde{a}(1-t)
\end{equation*}
since $\nu<1$. As a consequence, we can apply the analysis of the previous sections to Problem \eqref{eq26} and, going back to the original problem, we obtain the same imperfect bifurcations of Figure \ref{Fig6} if we represent on the vertical axis $u(1-\a)$ instead of $u(\a)$.

\section*{Acknowledgement}
The author wishes to thank Prof. J. L\'opez-G\'omez for the enlightening discussions and suggestions during the preparation of this work.

\end{document}